\documentclass[12pt]{amsart}

\usepackage{mathrsfs}
\usepackage{amsmath,amsthm,amssymb}
\font\teneufm=eufm10
\font\seveneufm=eufm7
\font\fiveeufm=eufm5
\newfam\frakturfam

\textfont\frakturfam=\teneufm
\scriptfont\frakturfam=\seveneufm
\scriptscriptfont\frakturfam=\fiveeufm

\def\bes{\begin{eqnarray*}}
\def\ees{\end{eqnarray*}}
\def\bee{\begin{eqnarray}}
\def\eee{\end{eqnarray}}
\def\la{\langle}
\def\ra{\rangle}

\def\a{\alpha}
\def\b{\beta}

\def\s{\sigma}
\def\f{\varphi}

\def\Z{\mathbf Z}
\def\0{\bar 0}
\def\1{\bar 1}

\def\H{{\mathcal H}}

\def\prf{{\it Proof.\ } }
\def\ctd{\hfill$\Box$}

\def\P{\mathcal P}

\newtheorem{thm}{Theorem}[section]
\newtheorem{cor}[thm]{Corollary}
\newtheorem{lem}[thm]{Lemma}
\newtheorem{prop}[thm]{Proposition}

\newtheorem{conj}[thm]{Conjecture}

\theoremstyle{definition}

\numberwithin{equation}{section}

\title{Representations of the Grassmann Poisson superalgebras}

\begin{document}
\date{}
\maketitle

\begin{center}
{\bf Ivan Shestakov}\footnote{Instituto de Matem\'atica e Estat\'\i stica, Universidade de S\~ao Paulo, 
Caixa Postal 66281, S\~ao Paulo - SP, \mbox{05315--970}, Brazil, 
and IMC SUSTech, Shenzhen, China,\\
e-mail:{\em shestak@ime.usp.br}}
and
{\bf Ualbai Umirbaev}\footnote{Department of Mathematics,
 Wayne State University,
Detroit, MI 48202, USA 
and Institute of Mathematics and Mathematical Modeling, Almaty, 050010, Kazakhstan,\\
e-mail: {\em umirbaev@wayne.edu}}
\end{center}

\begin{abstract} 
We prove that every irreducible Poisson supermodule over the Grassmann Poisson superalgebra $G_n$ over a field of characteristic different from $2$ is isomorphic to the regular Poisson supermodule $\mathrm{Reg}\,G_n$ or to its opposite supermodule.  Moreover,  every unital Poisson supermodule over $G_n$ is completely reducible.  If $P$ is a unital Poisson superalgebra which contains $G_n$ with the same unit then $P\cong Q\otimes G_n$ for some Poisson superalgebra $Q$.
Furthermore, we classify the supermodules over $G_n$ in the category of dot-bracket superalgebras with Jordan brackets, and we prove that every irreducible Jordan supermodule over the Kantor double  $\mathrm{Kan}\,G_n$ is isomorphic to the supermodule $\mathrm{Kan}\,V$, where $V$ is an irreducible dot-bracket supermodule with a Jordan bracket over $G_n$.
\end{abstract}

\noindent {\bf Mathematics Subject Classification (2020):} 17B63, 17B60, 17C70, 17A70

\noindent {\bf Key words:} Poisson superalgebra, Poisson supermodule,  Grassmann algebra, Jordan superalgebra,  Kantor double,  Jordan bracket, contact Lie bracket

\tableofcontents

\section{Introduction}

\hspace*{\parindent}
Let $G_n$ be the Grassmann algebra over a vector space of dimension $n$.  It  has a natural $\Z_2$-grading under which it forms a commutative superalgebra.  Moreover, it has also a super-anticommutative bracket ({\em a Poisson bracket}) and under the associative supercommutative multiplication ({\em a dot product} ) and this bracket it forms a {\em Poisson superalgebra}.  Over a field of characteristic zero every finite dimensional simple Poisson superalgebra is isomorphic to $G_n,\,n\geq 2$ \cite{Cheng}.

We first study representations of $G_n$ in  the category of Poisson superalgebras. 
It occurs that every irreducible Poisson supermodule over $G_n$ is isomorphic to the regular supermodule $\mathrm{Reg}\, G_n$ or to its parity-opposite module.  Moreover, every unital Poisson $G_n$-supermodule is completely reducible.  Using this facts, we prove the following Coordinatization Theorem for  $G_n$:

{\em Let $P$ be a Poisson superalgebra that contains $G_n$ with the same unit. 
Then there exists a Poisson subsuperalgebra $A$ of $P$ such that $P\cong A\otimes G_n$.}

This is an analogue of coordinatization theorems for different classes of algebras and superalgebras starting with the classical Wedderburn theorem for matrix algebras (see \cite{Jac1, Kap, LDS1,LDS2, LSSh,MSZ, MZ0,Mac,PS})

\smallskip

The superalgebra $G_n$ plays an important role in the theory of Jordan superalgebras: due to the {\em Kantor double process} with any superalgebra $G_n$ one can associate  a simple Jordan superalgebra $\mathrm{Kan}\, G_n$.  The Kantor construction $\mathrm{Kan}$ is functorial, and one can associate with any Poisson $G_n$-supermodule $V$  a Jordan supermodule $\mathrm{Kan}\,V$ over $\mathrm{Kan}\,G_n$ which is irreducible if $V$ is so. There was a conjecture that every irreducible Jordan supermodule over $\mathrm{Kan}\,G_n$ can be obtained  in this way.   It follows from our classification of irreducible Poisson supermodules over $G_n$ and from the results of  \cite{FSh,MZ} that it is not true: the irreducible Jordan $\mathrm{Kan}\,G_n$-supermodules form a  family parametrized by the scalars from the ground field $F$.
\smallskip

Fortunately,  the functor $\mathrm{Kan}$ can be applied not only to Poisson superalgebras but to any ``dot-bracket'' superalgebra $A$,   that is,  a superalgebra with an associative and commutative ``dot-multiplication''  $a\cdot b$ and a super-anticommutative bracket $\{a,b\}$. If the resulting commutative superalgebra $\mathrm{Kan}(A)$ is Jordan then the bracket $\{,\}$ is called a {\em Jordan bracket}. 

Thus  we decide to classify the supermodules over $G_n$ in the category of dot-bracket superalgebras with Jordan brackets.  It occurs that in this case every irreducible Jordan supermodules over the Kantor double  $\mathrm{Kan}\,G_n$ is isomorphic to the supermodule $\mathrm{Kan}\,V$ where $V$ is an irreducible dot-bracket supermodule with Jordan bracket over $G_n$.

It worth to be noticed that in fact we considered not Jordan brackets but so called {\em Lie contact brackets} which due to \cite{CK} are in one-to-one correspondence with Jordan brackets but are easier to deal with.

\section{Poisson superalgebras and supermodules}

\hspace*{\parindent}

We begin by reviewing some standard notions and facts needed for the proofs of the main results.  All (super)algebras and (super)modules are considered over a field $F$ of characteristic different from $2$.

A vector {\em superspace} $V=V_0\oplus V_1$ is a $\mathbb{Z}_2$-graded space.
If $v\in V_{\a}$, where $\a\in \mathbb{Z}_2=\{0,1\}$, we say that $\a$ is the {\em parity} of $v$ and denote it by $|v|$.

A vector superspace $P=P_0\oplus P_1$ over a field $F$ endowed with two bilinear
operations $x\cdot y$ (a multiplication) and $\{x,y\}$ (a Poisson bracket) is called
{\em a Poisson superalgebra}  if $P$ is a commutative
associative superalgebra under $x\cdot y$:
\bes
(x\cdot y)\cdot z=x\cdot (y\cdot z),\\
(x\cdot y)=(-1)^{|x||y|}(y\cdot x);
\ees
 $P$ is a Lie superalgebra under
$\{x,y\}$:

\bes
\{x,y\}=-(-1)^{|x||y|}\{y,x\}, \\
\{x,\{y,z\}\}=\{\{x,y\},z\}+(-1)^{|x||y|}\{y,\{x,z\}\};
\ees
 and $P$ satisfies the Leibniz rule:
\bes
\{x,
y\cdot z\}=\{x,y\}\cdot z + (-1)^{|x||y|} y\cdot \{x,z\}
\ees
for all $x,y,z\in P_0\cup P_1$. 

The Grassmann algebra $G=G_n$ is the associative algebra with identity $1$ generated by $e_1,\ldots,e_n$ and defined by the relations 
\bes
e_ie_j=-e_je_i,\,e_i^2=0 \ \  {\text for \  all}  \ \ 1\leq i\neq j\leq n.
\ees
  It has a basis 
\bes
1,\,e_{i_1}e_{i_2}\cdots e_{i_k},\,1\leq i_1<i_2<\cdots<i_k\leq n.  
\ees

If we set $|e_i|=1$ for all $i$,  then 
\bes
G =G_0\oplus G_1
\ees
 becomes a commutative and associative superalgebra, where $G_0$ and $G_1$ are the linear spans of all monomials of even and odd lengths, respectively. Moreover, it is a free superalgebra in the odd variables $e_1,e_2,\ldots,e_n$. The commutative superalgebra $F[x_1,\ldots,x_m]\otimes G_n$, where the polynomial algebra $F[x_1,\ldots,x_m]$ is regarded as a superalgebra with $|x_i|=0$ for all $i$, is a free commutative and associative superalgebra with even generators $x_1,\ldots,x_m$ and odd generators $e_1,e_2,\ldots,e_n$.

For Poisson superalgebras $P=P_0\oplus P_1$ and $Q=Q_0\oplus Q_1$ their tensor product $P\otimes Q$ is  defined as the vector superspace 
\bes
P\otimes Q = (P_0\otimes Q_0\oplus P_1\otimes Q_1)\oplus  (P_0\otimes Q_1\oplus P_1\otimes Q_0)
\ees
 with the following product and bracket
\bes
p\otimes q\cdot p_1\otimes q_1&=&(-1)^{|q||p_1|} pp_1\otimes qq_1,\\
\{p\otimes q, p_1\otimes q_1\}&=&(-1)^{|q||p_1|} (pp_1\otimes\{q,q_1\}+\{p,p_1\}\otimes qq_1).
\ees

Here are some important examples of Poisson (super)algebras.

$(1)$ {\em Symplectic Poisson algebra $P_m$}. For each $m$ the algebra $P_m$ is the polynomial algebra 
\bes
F[x_1,\ldots,x_m,y_1, \ldots,y_m]
\ees
 endowed with the Poisson bracket
\bes
\{f,g\}= \Sigma_{i=1}^m (\frac{\partial f}{\partial x_i}\frac{\partial g}{\partial y_i}-
\frac{\partial f}{\partial y_i}\frac{\partial g}{\partial x_i}).
\ees

$(2)$ {\em The Grassman Poisson superalgebra $G_n$} is the associative and commutative superalgebra $G_n$ endowed with the Poisson (super)bracket
\bes
\{f,g\}=(-1)^{|f|} \Sigma_{i=1}^n \frac{\partial f}{\partial e_i}\frac{\partial g}{\partial e_i}, 
\ees
where 
\bes
\tfrac{\partial }{\partial e_{i_s}}(e_{i_1}\cdots e_{i_s}\cdots e_{i_k})=(-1)^{s-1}e_{i_1}\cdots e_{i_{s-1}}e_{i_{s+1}}\cdots e_{i_k}.
\ees

$(3)$ {\em Poisson superalgebra $P_m\otimes G_n$}. Above described two brackets can be extended to the commutative superalgebra $P_m\otimes G_n$ by
\bes
\{f,g\}= \Sigma_{i=1}^m (\frac{\partial f}{\partial x_i}\frac{\partial g}{\partial y_i}-
\frac{\partial f}{\partial y_i}\frac{\partial g}{\partial x_i})+(-1)^{|f|} \Sigma_{i=1}^n \frac{\partial f}{\partial e_i}\frac{\partial g}{\partial e_i}.
\ees

$(4)$ {\em Symmetric Poisson algebra $PS(\mathfrak{g})$}. Let $\mathfrak{g}=\mathfrak{g}_0\oplus \mathfrak{g}_1$ be
a Lie superalgebra, $f_1,f_2,\ldots,f_k,\ldots$ be a linear basis of $\mathfrak{g}_0$, and  $g_1,g_2,\ldots,g_s,\ldots$ be a linear basis of $\mathfrak{g}_1$.
Then $PS(\mathfrak{g})$ is the commutative associative superalgebra
\bes
F[f_1,f_2,\ldots,f_k,\ldots]\otimes G(g_1,g_2,\ldots,g_s,\ldots), 
\ees
where $G(g_1,g_2,\ldots,g_s,\ldots)$ is the Grassmann algebra in the variables $g_1,g_2,\ldots,g_s,\ldots$, 
 with the Poisson bracket determined by 
 \bes \{x,y\}=[x,y]
\ees
 for all $x,y\in \{f_1,f_2,\ldots,f_k,\ldots,g_1,g_2,\ldots,g_s,\ldots\}$, where $[x,y]$
is the multiplication of the Lie superalgebra $\mathfrak{g}$.

It is well known that the symplectic Poisson algebra $P_m$ is simple. The following result is also well known.

\begin{prop}\label{prop2}\cite{Kantor}
The Grassmann Poisson superalgebra $G_n$ for $n>1$  is a simple Poisson superalgebra.
\end{prop}
\begin{cor}
 The Poisson superalgebra  $P_m\otimes G_n$ for $n>1$ is simple. 
\end{cor}

Every simple Lie superalgebra, regarded as a Poisson superalgebra with trivial multiplication, is a simple Poisson algebra. Every finite dimensional simple Poisson superalgebra $P$ with an identity ( or with a nontrivial multiplication) over an algebraically closed field of characteristic zero is isomorphic to $G_n$.  This follows from the fact that $P$ is simple as a Poisson superalgebra if and only if, when regarded as a commutative superalgebra, it is differentially simple with respect to the derivations $h_a:x\mapsto \{a,x\}$, where $a\in P$; and any differentially simple commutative superalgebra over an algebraically closed field of characteristic zero is isomorphic to $G_n$ by \cite[theorem 4.1]{Cheng}.

A vector superspace $V=V_0\oplus V_1$ is called a {\em Poisson supermodule} over a Poisson superalgebra $P=P_0\oplus P_1$ if  the two even linear mappings are defined
\bes
m, h : P\rightarrow \mathrm{End}\, V,
\ees
which define the two actions of $P$ on $V$:
\bes
v\cdot a=v\, m(a),\ \{v,a\}=v\, h(a),
\ees
such that the {\em split null extension} $E(P,V)=P\oplus V$
with the operations 
\bes
(a+v)(b+u)&=&ab+(v\cdot b+(-1)^{|u||| a|}u\cdot a),\\
\{a+v,b+u\}&=&\{a,b\}+(\{v, b\}+(-1)^{|u||a|}\{u,a\})
\ees
becomes a Poisson superalgebra with the grading
\bes
E(P,V)_0=P_0\oplus V_0, \ E(P,V)_1=P_1\oplus V_1.
\ees

It is easy to see that the mappings $m,h$ define a Poisson supermodule structure on $V$ if and only if they satisfy the following identities:
\bee
m(a\cdot b)&=&m(a)m(b), \label{id1}\\
m(\{a,b\})&=&m(a)h(b)-(-1)^{|b||a|}h(b)m(a),\label{id2}\\
h(a\cdot b)&=&h(a)m(b)+(-1)^{|b||a|}h(b)m(a),\label{id3}\\
h(\{a,b\})&=&h(a)h(b)-(-1)^{|b|||a|}h(b)h(a).\label{id4}
\eee
In this case  the pair $(m,h)$ is called a {\em representation} of the superalgebra $P$ on the module $V$.  Clearly,  the notions of module and representation mutually define each other. 

In a standard way (see, for instance \cite{Jacobson, KSU, U96,U12}) it is proved that there exists the universal 
associative superalgebra $U(P)$ ({\em the universal multiplicative envelope of $P$}) and the linear mappings $\mathcal{M},\mathcal{H}:\P\rightarrow U(P)$ that satisfy the above identities for $m,h$ and such that for any representation $(m,h):P\rightarrow \mathrm{End}\,V$
there exists a unique homomorphism $\phi :U(P)\rightarrow \mathrm{End}\,V$ satisfying the equalities
\bes
m=\phi\circ \mathcal M,\ h=\phi\circ \H.
\ees

The category of Poisson $P$-supermodules is isomorphic to the category of associative right $U(P)$-supermodules. 

Every Poisson superalgebra $P$ is itself a supermodule over $P$. This module is denoted by $\mathrm{Reg}\,P$, and the corresponding representation is called the {\em regular representation} of $P$. 

For any Poisson supermodule $V=V_0\oplus V_1$ over  $P=P_0\oplus P_1$ the {\em opposite} $P$-supermodule 
\bes
V^{\mathrm{op}}=V^{\mathrm{op}}_0\oplus V^{\mathrm{op}}_1, 
\ees
where $V^{\mathrm{op}}_0=V_1$ and $V^{\mathrm{op}}_1=V_0$, 
 is defined by 
\bes
v^{\mathrm{op}}\cdot p=(v\cdot p)^{\mathrm{op}}, \{v^{\mathrm{op}}, p\}=\{v, p\}^{\mathrm{op}}
\ees 
for any $v\in V_0\cup V_1$ and $p\in P_0\cup P_1$. The identity map 
\bes
\mathrm{Id} : V\to V^{\mathrm{op}} (v_0+v_1\mapsto v_1^{op}+v_0^{op})
\ees 
is an odd isomorphism of $P$-modules. In general, there is no even isomorphism between $V$ and $V^{\mathrm{op}}$.

\section{Poisson representations of $G_n$}

\hspace*{\parindent}

In this section we describe the structure of the universal enveloping algebra $U(G_n)$ and describe all finite dimensional representations of $G_n$. 

\begin{thm}\label{thm1}
1. The multiplicative enveloping superalgebra $U(G_n)$ is isomorphic to the Clifford superalgebra $Cl(W)$ of an odd vector space  $W=W_1$ of dimension $2n$. 

2.
Every irreducible unital Poisson $G_n$-supermodule is isomorphic to the regular supermodule $\mathrm{Reg}\,G_n$  or to its opposite supermodule.   

3.  Any unital Poisson module over $G_n$ is completely reducible and is isomorphic to a direct sum of modules $\mathrm{Reg}\,G_n$ and $(\mathrm{Reg}\,G_n)^{op}$.
\end{thm}
\prf
Consider in $U(G_n)$ the subspace $W$ spanned by the odd elements 
$$
v_1=\mathcal M(e_1),\ldots,v_n=\mathcal M(e_n); v_{n+1}=\H(e_1),\ldots,v_{2n}=\H(e_n).
$$
It follows from the identities \eqref{id1} - \eqref{id4} that the space $W$ generates the algebra $U(P)$.  Moreover, we have 
\bes
v_iv_j+v_jv_i&=&\mathcal M(e_i)\mathcal M(e_j)+\mathcal M(e_j)\mathcal M(e_i)=\mathcal M(e_ie_j+e_je_i)=0,\\
v_{n+i}v_{n+j}+v_{n+j}v_{n+i}&=&\H(e_i)\H(e_j)+\H(e_j)\H(e_i)=\H\{e_i,e_j\}=\H(-\delta_{ij}\cdot  1)=0, \\
v_iv_{n+j}+v_{n+j}v_i&=&\mathcal M(e_i)\H(e_j)+H(e_j)\mathcal M(e_i)
=\mathcal M(\{e_i,e_j\}=\mathcal M(-\delta_{ij}\cdot 1)=-\delta_{ij}, 
\ees
for all $i,j\leq n$. 

Define on the space $W$ the symmetric bilinear form $f(x,y)$ as follows: 
\bes
f(v_i,v_j)&=&0 \hbox{ if } i,j\leq n \hbox{ or } i,j>n;\\
f(v_i,v_{n+j})&=&f(v_{n+j},v_i)=-\delta_{ij}.
\ees
Clearly, the form $f(x,y)$ is nondegenerated on $W$.  Moreover, the above relations show that for any $u,w\in W$ we have
\bes
uw+wu=f(u,w)\cdot 1.
\ees
This proves that the algebra $U(P)$ is isomorphic to the Clifford algebra $\mathrm{Cl}(W,f)$ of the form $f$ on the space $W$.  

\smallskip

The algebra $\mathrm{Cl}(W,f)$ has a basis 
$$
1,\, v_{i_1}v_{i_2}\cdots v_{i_k},\ 1\leq i_1<i_2<\cdots<i_k\leq 2n.
$$
It has a $\Z_2$-grading determined by the odd subspace $W$:  the even part $\mathrm{Cl}(W,f)_0$ is spanned by 1 and the products of even length, and the odd part $\mathrm{Cl}(W,f)_1$ is spanned by the products of odd length.

Since $\dim W=2n$ is even, the algebra $\mathrm{Cl}(W,f)$ is simple and is isomorphic to the matrix algebra $M_{2^{n}}(F)$.   As a superalgebra, it is isomorphic to the superalgebra $M_{2^{n-1},2^{n-1}}(F)$.  
It is easy to see that, up to changing of parity,  it has only one irreducible  supermodule.  Clearly, the regular module $\mathrm{Reg}\,G_n$ and its opposite are irreducible.  Consequently, every irreducible Poisson module over $G_n$ is isomorphic to $\mathrm{Reg}\,G_n$ or to $(\mathrm{Reg}\,G_n)^{op}$. 

To prove the last statement of the theorem, it suffices to notice that  any  graded module over $M_{2^{n}}(F)$ is completely reducible as a module, and  all its irreducible components are isomorphic to $\mathrm{Reg}\,G_n$. 
\ctd

\section{Coordinatization theorem}

\hspace*{\parindent}

Let $I_n=\{1,2.\ldots,n\}$. Let  $I\subseteq I_n$. If $I=\{i_1,\cdots,i_k\,|1\leq i_1<\cdots<i_k\leq n\}$ then set $e_I=e_{i_1}\cdots e_{i_k}$. In particular, $e_{\emptyset}=1$. The set of all such elements 
\bee\label{f1}
e_I, \ \ I\subseteq I_n, 
\eee
 is a linear basis of $G_n$.

We have 
\bes
G_n\otimes G_m\cong G_{n+m} 
\ees
for any $m,n\geq 0$. This is an analogue of the well known isomorphism 
\bes
M_n(F)\otimes M_m(F)\cong M_{nm}(F). 
\ees
It is well known if $A$ is finite dimensional associative algebra containing $M_n(F)$ with the same unit then there exists a subalgebra $B$ of $A$ such that $A\cong B\otimes M_n(F)$. We prove an analogue of this result for $G_n$ in the case of Poisson superalgebras. 

\begin{thm}\label{thm2}
Let $P$ be a Poisson superalgebra that contains $G_n$ with the same unit. 
Then there exists a Poisson subsuperalgebra $A$ of $P$ such that $P\cong A\otimes G_n$.
\end{thm}
\prf
Let $A=\{a\in P\,|\,\{a,g\}=0 \hbox{ for any } g\in G_n\}$.  It follows from the Leibniz and super-Lie identities that $A$ is a subsuperalgebra of $P$.  

Consider $P$ as a Poisson $G_n$-module. By theorem \ref{thm1}, it is completely reducible and $P=\oplus_i P_i$,  where $P_i\cong \mathrm{Reg}\,G_n$ or $P_i\cong (\mathrm{Reg}\,G_n)^{op}$ for  all $i$.  Let  $a_i\in P_i$ be the generator of $P_i$ that corresponds to $1\in \mathrm{Reg}\,G_n$ or to $1^{op}\in (\mathrm{Reg}\,G_n)^{op}$.  Clearly, all $a_i\in A$,  which proves that $P\subseteq A\cdot G_n$.

To prove the isomorphism $A\cdot G_n\cong A\otimes G_n$ of the vector spaces, we need to prove that the elements of the basis \eqref{f1} of $G_n$ are linearly independent over $A$.  Assume that 
\bee\label{id5}
\sum_{I\subseteq I_n} a_I\cdot e_I=0. 
\eee
Choose a basis element $e_I$ with nonzero coefficient $a_I$ and having the minimal number of factors $e_i$. Then every other element $e_J$ contains a factor $e_j$ such that $j\not\in I$.  Multiplying  the relation \eqref{id5} successively by such elements $e_{j_1},e_{j_2},\ldots,e_{j_m}$, we eventually get $a_I\cdot e_Ie_{j_1}\cdots e_{j_m}=0$. 
 
 \smallskip
 Furthermore, we have
 \bes
 0&=&\{a_Ie_Ie_{j_1}\cdots e_{j_m},e_{j_m}\}\\
 &=&-\{a_Ie_Ie_{j_1}\cdots e_{j_{m-1}},e_{j_m}\}\cdot e_{j_m}+a_Ie_Ie_{j_1}\cdots e_{j_{m-1}}\cdot\{e_{j_m},e_{j_m}\}\\
 &=&-\{a_Ie_Ie_{j_1}\cdots e_{j_{m-1}},e_{j_m}\}\cdot e_{j_m}-a_Ie_Ie_{j_1}\cdots e_{j_{m-1}}\\
 &=&\cdots\\
 &=&\pm \{a_I,e_{j_{m}}\}e_Ie_{j_1}\cdots e_{j_{m-1}}e_{j_m}-a_Ie_Ie_{j_1}\cdots e_{j_{m-1}}\\
 &=&-a_Ie_Ie_{j_1}\cdots e_{j_{m-1}}.
\ees
Continuing in this way, we eventually get $a_I=0$.

\smallskip

Finally, we have 
\bes
(a\cdot g)(b\cdot h)&=&(-1)^{|g|| b|} ab\cdot gh,\\
\{a\cdot g,b\cdot h\}&=&\{a\cdot g,b\}\cdot h+(-1)^{|b||h|}\{a\cdot g,h\}\cdot b\\
&=&(-1)^{|g||b|}\{a,b\}\cdot gh+(-1)^{|b||g|}ab\cdot \{g,h\}\\
&=&(-1)^{|g||b|}(\{a,b\}\cdot gh+ab\cdot \{g,h\}),
\ees
for all $a,b\in A,\ g,h\in G_n$. This proves the isomorphism $P\cong A\otimes G_n$ of Poisson superalgebras. 
\ctd

\section{The Kantor double, Jordan brackets, and contact Lie brackets}.

\hspace{\parindent}

I. Kantor  \cite{Kantor} introduced a functor from the category of Poisson (super)algebras to the category of Jordan superalgebras.  Let $P=P_0\oplus P_1$ be a Poisson superalgebra with multiplication $ab$ and bracket $\{a,b\}$, and let $\bar P$ be an isomorphic copy of the vector superspace $P$.  Consider the vector space direct sum 
\bes
\mathrm{Kan}(P)=P\oplus \bar P
\ees
 and define a multiplication $\cdot$ on it by setting 
\bes
a\cdot b&=&ab,\\
a\cdot \bar b&=&\overline{ab},\\
\bar a\cdot b&=&(-1)^{|b|}\overline{ab},\\ 
\bar a\cdot\bar b&=&(-1)^{|b|}\{a,b\}, 
\ees
for all $a,b\in P_0\cup P_1$. 
Define a grading on $\mathrm{Kan}(P)$ by setting 
$$
\mathrm{Kan}(P)_0=P_0\oplus \bar P_1,\ \mathrm{Kan}(P)_1=P_1\oplus \bar P_0. 
$$
Then $\mathrm{Kan}(P)$  becomes a Jordan superalgebra (see \cite{Kantor}).  

The mapping $P\mapsto \mathrm{Kan}(P)$ is functorial; in particular, if $P$ is a simple Poisson superalgebra then $\mathrm{Kan}(P)$ is a simple Jordan superalgebra.  The functor $\mathrm{Kan}$ can be extended to the associated categories of modules:
\bes
\mathrm{Kan}: P\hbox{-$_{\mathrm{Pois}}$mod}\rightarrow \mathrm{Kan}(P)\hbox{-$_{\mathrm{Jord}}$mod},
\ees
constructing for a Poisson $P$-(super)module $V$ a Jordan (super)module $\mathrm{Kan}(V)$.

A conjecture made by Efim Zelmanov and the first author states that every irreducible Jordan supermodule over $\mathrm{Kan}(P)$ is of the form $\mathrm{Kan}(V)$ for some irreducible Poisson $P$-supermodule   $V$. Theorem \ref{thm1} provides a negative answer to this conjecture.  In fact, in \cite{FSh, MZ} irreducible Jordan supermodules over $\mathrm{Kan}(G_n)$ were constructed that are not isomorphic to $\mathrm{Kan}(\mathrm{Reg}\,G_n)$.

\medskip
Recall that the functor $\mathrm{Kan}$ can be applied not only to Poisson superalgebras but to any ``dot-bracket'' superalgebra $A$,   that is,  a superalgebra with an associative and commutative ``dot-multiplication''  $a\cdot b$ and a super-anticommutative bracket $\{a,b\}$. If the resulting commutative superalgebra $\mathrm{Kan}(A)$ is Jordan then the bracket $\{,\}$ is called a {\em Jordan bracket}. 

D. King and K. McCrimmon  proved  \cite{KM1,KM2} that a bracket $\{a,b\}$ is a Jordan bracket if and only if it satisfies the identities
\bes
\{a,bc\}& =& \{a,b\}c + (-1)^{|a||b|} b\{a,c\} - D(a)bc,\\
J(a,b,c) &:=& \{\{a,b\},c\} + (-1)^{|a||b| + |a||c| } \{\{b,c\},a\}
+(-1)^{|a||c| + |b||c| } \{\{c,a\},b\}  \nonumber \\
& =& -\{a,b\}D(c) - (-1)^{|a||b| +|a||c| } \{b,c\}D(a) - (-1)^{|a||c| + |b||c| } \{c,a\}D(b),\\
\{\{x,x\},x\}&=&-\{x,x\}D(x), 
\ees
where $x$ is odd and $D(a)=\{a,1\}$.
The last identity  is needed only in characteristic 3 case,  otherwise it follows
from the previous one.  If $D=0$ we get a Poisson bracket.

N. Cantarini and V. Kac \cite{CK} showed that all the Kantor  doubles $\mathrm{Kan}(A)$,
which are Jordan superalgebras,  can be obtained from a contact Lie bracket on the superalgebra $A$.  By definition, a {\em contact Lie bracket} is a Lie superalgebra bracket \{·,·\} satisfying the generalized Leibniz rule
\bes
\{a,bc\} = \{a,b\}c + (-1)^{|a||b|} b\{a,c\} + D(a)bc, 
\ees
where $D(a)=\{1,a\}$ is an even derivation of the product and the bracket. 

For a contact Lie bracket $\{,\}$,   the new bracket 
\bee\label{Jord}
\la a,b\ra=\{a,b\}-\tfrac12(a \{1,b\}-\{1,a\}b)
\eee
is a Jordan bracket.  Conversly,  for a Jordan bracket $\la,\ra$,  the new bracket 
\bee\label{cont}
\{a,b\}=\la a,b\ra+ (a\la 1,b\ra-\la 1,a\ra b)
\eee
is a contact Lie bracket.  

\smallskip
It is easy to see that  any finite dimensional unital associative commutative superalgebra $A$ over an algebraically closed field $F$ of zero characteristic with a Jordan or contact Lie bracket which is simple as a dot-bracket algebra is differentially simple and hence by \cite{Cheng} is isomorphic to the algebra $G_n$ with the above defined Poisson bracket. 

But the structure of the universal multiplicative enveloping algebra $U(G_n)$ and of irreducible supermodules over $G_n$ depends on the category in which this algebra is considered.

We are going to describe the structure of irreducible supermodules over $G_n$ in the category of superalgebras with a Jordan brackets.  In view of the above equivalence between Jordan and contact Lie brackets, we prefer to work first with contact Lie brackets, as they are easier to handle.

\section{Representations of $G_n$ as a superalgebra with a contact Lie bracket}.

\hspace{\parindent}

A vector superspace $V=V_0\oplus V_1$ is  a {\em  supermodule with a contact Lie bracket}  over a superalgebra $P=P_0\oplus P_1$ with a contact Lie  bracket if  the two 
actions $m$ and $h$ of $P$ on $V$ satisfy identities \eqref{id1}, \eqref{id4}, and  the identities
\bee
m(\{a,b\})&=&m(a)h(b)-(-1)^{|a||b|}h(b)m(a)+m(a)m(\{1,b\}),\label{id6}\\
h(ab)&=&h(a)m(b)+(-1)^{|a||b|}h(b)m(a)-h(1)m(a)m(b).\label{id7}
\eee

Since $G_n$ is a Poisson algebra, we have $\{1,b\}=0$ for all $b\in G_n$. Consequently,   for contact Lie supermodules over $G_n$ the identity \eqref{id6} 
coinsides with \eqref{id2}. Thus, in this section we can use the identities \eqref{id1},  \eqref{id2},\eqref{id4}, and \eqref{id7}.

Notice also that $\H(1)$ lies in the center of the universal multiplicative enveloping algebra $U_{CLie}(G_n)$ in the category of superalgebras with contact Lie brackets.  In fact,  using \eqref{id2} and \eqref{id4}, we get 
\bes
[\mathcal M(a),\H(1)]&=&\mathcal M(\{a,1\})=0,\\ \
[\H(a),\H(1)]&=&\H_{\{a,1\}}=0.
\ees

We are going all  irreducible finite-dimensional $G_n$-supermodules in the category of contact Lie brackets. First give some examples. 

Let $\overline{G_n}$ be an isomorphic copy of $G_n$.  Define on the vector space $\overline{ G_n}$ two actions of $G_n$: 
\bes
\bar e_I\cdot e_J&=&\overline{e_I\cdot e_J},\\ 
\{\bar e_I,e_J\}&=&\overline{\{e_I,e_J\}}+\b\,(|J|-2) \,\overline{e_I\cdot e_J}, 
\ees
where $I,J\subseteq I_n$ and $|J|$ is the number of elements of $J$. 

\begin{prop} \label{prop3}
The space $\overline{G_n}$ is a supermodule for the dot-bracket superalgebra $G_n$ in the category of superalgebras with contact Lie brackets.
\end{prop}
\prf  We have to check the identities \eqref{id1},  \eqref{id2}, \eqref{id4}, and \eqref{id7}. Let $I,J,K\subseteq I_n$.  We have  
\bes
(\bar e_I \cdot e_J)\cdot e_K=\bar (e_I \cdot e_J)\cdot e_K=\bar (e_I \cdot e_J\cdot e_K)= (\bar e_I) \cdot (e_J\cdot e_K), 
\ees
whih proves \eqref{id1}. 

Furthermore, 
\bes
\{\bar e_I\cdot e_J,e_K\}&=&\{\overline{e_I\cdot e_J},e_K\}\\
&=&\overline{\{e_I\cdot e_J,e_K\}}+\b(|K|-2)\overline{e_I\cdot e_J\cdot e_K}\\
&=&\overline{e_I\cdot\{e_J,e_K\}}+(-1)^{|J||K|}\overline{\{e_I,e_K\}\cdot e_J}
+\b(|K|-2)\overline{e_I\cdot e_J\cdot e_K}\\
&=&\bar e_I\cdot \{e_J,e_K\}+(-1)^{|J||K|}\overline{\{e_I,e_K\}}\cdot e_J+\b(|K|-2)\overline{e_I\cdot e_J\cdot e_K}\\
&=&\bar e_I\cdot \{e_J,e_K\}+(-1)^{|J||K|}(\{\bar e_I,e_K\}-\b\,(|K|-2)\,\overline{e_I\cdot e_K})\cdot e_J\\
&+&\b\,(|K|-2)\,\overline{e_I\cdot e_J\cdot e_K}\\
&=&\bar e_I\cdot \{e_J,e_K\}+(-1)^{|J||K|}\{\bar e_I,e_K\}\cdot e_J.
\ees
This proves \eqref{id2}.

Direct calculations give that 
\bes
\{\{\bar e_I,e_J\},e_K\}&=&\{\overline{\{e_I,e_J\}}+\b (|J|-2)\,\overline{e_I\cdot e_J},e_K\}\\
&=&\overline{\{\{e_I,e_J\},e_K\}}+\b(|K|-2)\,\overline{\{e_I,e_J\}\cdot e_K}\\
&+&\b(|J|-2)(\overline{\{e_I\cdot e_J,e_K\}}+\b (|K|-2)\,\overline{e_I\cdot e_J\cdot e_K})\\
&=&\overline{\{\{e_I,e_J\},e_K\}}+\b(|K|-2)\,\overline{\{e_I,e_J\}\cdot e_K}\\
&+&\b(|J|-2)\,(\overline{e_I\cdot \{e_J,e_K\}}+(-1)^{|J||K|}\overline{\{e_I,e_K\}\cdot e_J}+\b (|K|-2)\,\overline{e_I\cdot e_J\cdot e_K}).
\ees
Symmetrically, 
\bes 
\{\{\bar e_I,e_K\},e_J\}&=&\overline{\{\{e_I,e_K\},e_J\}}+\b(|J|-2)\,\overline{\{e_I,e_K\}\cdot e_J}\\
&+&\b(|K|-2)(\overline{e_I\cdot \{e_K,e_J\}}+(-1)^{|J||K|}\overline{\{e_I,e_J\}\cdot e_K}+\b (|J|-2)\,\overline{e_I\cdot e_K\cdot e_J}). 
\ees
It is clear that $\{e_J,e_K\}$ is a linear combination of monomials of length $|J|+|K|-2$. Then 
\bes
\{\bar e_I,\{e_J,e_K\}\}=\overline{\{e_I,\{e_J,e_K\}\}}+\b(|J|+|K|-4)\,\overline{e_I\cdot\{e_J,e_K\}}. 
\ees
Using these expressions, we get 
\bes
\{\{\bar e_I,e_J\},e_K\}-(-1)^{|J||K|}\{\{\bar e_I,e_K\},e_J\}-\{\bar e_I,\{e_J,e_K\}\}=0,
\ees
which proves \eqref{id4}.

Notice that $\{\bar e_I,1\}=-2\b e_I$. Then \eqref{id7} can written as 
\bee\label{id8}
\{\bar e_I,e_J\cdot e_K\}-\{\bar e_I,e_J\}\cdot e_K-(-1)^{|J||K|}\{\bar e_I,e_K\}\cdot e_J
-2\b \,\bar e_I\cdot (e_J\cdot e_K)=0.
\eee
We have
\bes
\{\bar e_I,e_J\}\cdot e_K&=&(\overline{\{e_I,e_J\}}+\b(|J|-2)\,\overline{e_I\cdot e_J})\cdot e_K\\
&=&\overline{\{e_I,e_J\}\cdot e_K}+\b \,(|J|-2)\,\overline{e_I\cdot e_J\cdot e_K}
\ees
and, symmetrically, 
\bes
\{\bar e_I,e_K\}\cdot e_J&=&\overline{\{e_I,e_K\}\cdot e_J}+\b \,(|K|-2)\,\overline{e_I\cdot e_K\cdot e_J}. 
\ees

Therefore,
\bes
\{\bar e_I,e_J\}\cdot e_K+(-1)^{|J||K|}\{\bar e_I,e_K\}\cdot e_J\\
=\overline{\{e_I,e_J\cdot e_K\}}+\b(|J|+|K|-4)\overline{e_I\cdot e_J\cdot e_K}\\
=\{\bar e_i,e_J\cdot e_K\}-2\b\, \overline{e_I\cdot e_J\cdot e_K}.
\ees
which proves \eqref{id8}.
\ctd

Notice that the structure of the module $\overline{G_n}$ depends on the choice of the scalar $\b\in F$.  For instance, for $\b=0$ we have $\overline{G_n}\cong \mathrm{Reg}\, G_n$.  We denote this module for a fixed $\b\in F$ as $G_n(\b)$.

In what follows, $V$ denotes an irreducible finite-dimensional $G_n$-supermodule in the category of contact Lie brackets.

\begin{lem}\label{lem1}
(1). $m(1)=Id_V,$

(2).  $h(1)=\alpha \cdot Id_V $ for some $\alpha\in F$.
\end{lem}
\prf
Let $V'=\{v-v\cdot 1\,|\,v\in V\}$.  By \eqref{id1} and \eqref{id2},  we have 
\bes
(v-v\cdot 1)\cdot a&=&v\cdot a-(v\cdot a)\cdot 1=v\cdot a-(v\cdot a)\cdot 1,\\
\{v-v\cdot 1,a\}&=&\{v,a\}-v\cdot\{1,a\}-\{v,a\}\cdot 1- v\cdot 1 \{1,a\}\\
&=&\{v,a\}-\{v,a\}\cdot 1
\ees
for any $a\in G_n$. This proves that $V'$ is a subsupermodule of $V$.  

Assume that $V'=V$. Then $v\cdot 1=0$ and $v\cdot a=v\cdot (1\cdot a)=0$ for any $v\in V,\,a\in G_n$. Hence $V\cdot G_n=(0)$.  Moreover, we have 
\bes
\{v,a\}=\{v,a\cdot 1\}=\{v,a\}\cdot 1+\{v,1\}\cdot a+\{1,v\}\cdot a=0,
\ees
and hence $\{V,G_n\}=(0)$, a contradiction.  Therefore, $V'=(0)$, completing the proof of (1).

\smallskip
The second statment follows from the fact that  $h(1)$ lies in the centralizer $Z_{G_n}(V)$ of the $G_n$-module $V$, which is a division algebra.  Since the field $F$ is algebraically closed,  we have $Z_{G_n}(V)=F$, hence 
$h(1)=\alpha\in F$.

\ctd

Notice that there exists a nonzero element $v\in V$ such that $v\cdot e_i=0$  for all  $i\in I_n=\{1,2,\ldots,n\}$.
 Indeed, let $v_n=v\cdot e_{I_n}$. If $v_n\neq 0$, then it satisfies the required property.
 If $v_n=0$,  let $I\subsetneqq I_n$ be a maximal subset such that $v\cdot e_I\neq 0$ but $v\cdot e_J=0$ for every $J\supsetneqq I.$ Then $v\cdot e_I$ satisfies the required property.
 
Let $v\in V$ be the element obtained above.  For any subset $I\subseteq I_n$ set 
\bes
v_I=\{\{\cdots\{v,e_{i_1}\},\cdots\},e_{i_k}\}
\ees
if $I=\{i_1<i_2<\ldots<i_k\}$. If $j\in I_n\setminus I$,  then denote by $s(I,j)$ the number of inversions in the sequence $i_1,\ldots,i_k,j$. 

We also fix $\a\in F$ satisfying the condition (2) of Lemma \ref{lem1}. 

\begin{lem}\label{lem2} Let $I=\{i_1<\ldots <i_k\}$. Then the following statements hold: 

(1) $v_I\cdot e_j=0$ for any $j\not\in I$,

(2) $v_I\cdot e_j=(-1)^{k-s-1}v_{I\setminus\{j\}}$ for  $j=i_s\in I$,

(3) $\{v_I,e_j\}=(-1)^{s(I,j)}v_{I\cup\{j\}}$ for $j\not\in I$,

(4) $\{v_I,e_j\}=(-1)^{k-s-1}\tfrac{\alpha}{2} v_{I\setminus \{j\}}$ for $j=i_s\in I$. 
\end{lem}
\prf
We proceed by induction on $k$.  We have 
\bes
v\cdot e_i&=&0,\  \{v,e_j\}=v_{\{j\}}
\ees 
by the definitions. Therefore, statements (1) and (3) of the lemma are hold for $k=0$. 
By \eqref{id2} and \eqref{id4}, we get 
\bes
\{v,e_i\}\cdot e_i&=&-\{v\cdot e_i,e_i\}+v\cdot\{e_i,e_i\}=-v,\\
2\{\{v,e_i\},e_i\}&=&\{v,\{e_i,e_i\}\}=-\{v,1\}=-\alpha v. 
\ees
This proves (2) and (4) for $k=1$.

Assume that the lemma is proved for $k'<k$.
Let $j=i_s\in I$. Then $v_I=(-1)^{k-s} \{v_{I'},e_j\}$, where $I'=I\setminus \{j\}$, by the induction proposition. 

By \eqref{id2} and \eqref{id4}, we get 
\bes
\{v_{I'},e_{j}\}\cdot e_j&=&-\{v_{I'}\cdot e_j,e_{j}\}+v_{I'}\cdot \{e_j,e_j\}=-v_{I'},\\
2\{\{v_{I'},e_{j}\},e_j\}&=&\{v_{I'},\{e_j,e_j\}\}=-\{v_{I'},1\}=-\alpha v_{I'}, 
\ees
and, consequently, 
\bes
v_I\cdot e_j=(-1)^{k-s-1}v_{I'}, \ \ \{v_I,e_j\}=(-1)^{k-s-1}\tfrac{\alpha}{2}v_{I'}. 
\ees
This proves the statements (2) and (4).  

Let $j\notin I$. By the induction proposition and \eqref{id2}, we get 
 \bes
v_I\cdot e_j&=&\{v_{I\setminus\{i_k\}},e_{i_k}\}\cdot e_j\\
&=&-\{v_{I\setminus\{i_k\}}\cdot e_j,e_{i_k}\}+v_{I\setminus\{i_k\}}\cdot \{e_j,e_{i_k}\}=0. 
\ees
This proves the statement (1). If $j>i_k$ then $\{v_I,e_j\}=v_{I\cup \{j\}}$ by the definition. If $j<i_k$, then by the induction proposition and \eqref{id4}, we get 
\bes
\{v_I,e_j\}&=&\{\{v_{I\setminus\{i_k\}},e_{i_k}\},e_j\}=-\{\{v_{I\setminus\{i_k\}},e_j\},e_{i_k}\}\\
&=&-(-1)^{s((I\setminus\{i_k\})\cup\{j\})}\{v_{(I\setminus\{i_k\})\cup \{j\}},e_{i_k}\}=(-1)^{s(I,j)}v_{I\cup \{j\}},
\ees
which proves the statement (4).   
\ctd

\begin{cor}\label{cor2} 
The module $V$ is spanned by the elements $v_I,\ I\subseteq I_n$.
\end {cor}
\prf The identities \eqref{id1},  \eqref{id2}, \eqref{id4}, and \eqref{id7} imply that the universal enveloping algebra $U_{CLie}(G_n)$ of  $G_n$ as a contact Lie bracket is generated by the elements $\H(1), \,\H(e_i),\,\mathcal M(e_i)$.  By Lemmas \ref{lem1} and \ref{lem2}, the space of elements spanned by all $v_I,\ I\subseteq I_n$, is closed under the action of $G_n$. Since $V$ is irreducible, 
it coincides with the span of $v_I,\ I\subseteq I_n$. 
\ctd
\begin{lem}\label{lem3}
The elements $v_I,  \ I\subseteq I_n$, are linearly independent.
\end{lem}
\prf
Assume that 
\bes
\sum_{I\subseteq I_n} \a_I v_I=0.
\ees
 Choose and fix a subset $I$ in this sum with nonzero coefficient $\a_I$ and with the maximal number of elements.  Applying the same discussions as in the proof of Theorem \ref{thm2} concerning the identity \eqref{id5}, and applying Lemma \ref{lem2}, we get that $\a_I=0$. 
\ctd

By passing to the opposite supermodule $V^{op}$,  if necessary, we may assume that when  $n$ is even the element $v$ is even, and when $n$ is odd the element $v$ is odd.

\begin{lem}\label{lem4} Suppose that $v\in V_0$ when $n$ is even and $v\in V_1$ when $n$ is odd. Then the map 
\bes
\f :  V\rightarrow G_n(\b),
\ees
for $\b=\tfrac{\a}{2}$, defined by 
\bes
\f: v_I\rightarrow (-1)^{k(n-1)+i_1+\ldots i_k}\bar e_{I'},
\ees
where $I=\{i_1<\ldots <i_k\}\subseteq I_n$ and $I'=I_n\setminus I$, is an isomorphism of $G_n$-supermodules. 
\end{lem}
\prf It is clear that $\f$ is an even isomorphism of vector superspaces. 
Set $w=\bar e _{I_n}$. Then $w\cdot e_j=0$ for all $j$. For any $I=\{i_1<\ldots <i_k\}\subseteq I_n$ set also 
\bes
w_I=\{\{\cdots\{w,e_{i_1}\},\cdots\},e_{i_k}\}. 
\ees
We have  
\bes
\{w, e_j\}= (-1)^{(n-1)+j}\bar e_{I_n\setminus \{j\}}=(-1)^{(n-1)+j}\bar e_{\{j\}'}. 
\ees
Continuing these calculations, we can get 
\bes
w_I=(-1)^{k(n-1)+i_1+\ldots i_k} \bar e_{I'}. 
\ees
Then $\f$ is defined by $\f(v_I)=w_I$ for all $I\subseteq I_n$. The statements of Lemma \ref{lem2} hold for all $w_I$, since $w$ satisfies the same conditons as $v$. This means that $\f$ preserves the actions. 
\ctd

This lemma implies the following theorem. 
\begin{thm}\label{thm3} 
Every irreducible finite dimensional contact Lie supermodule over the Grassmann Poisson algebra $G_n$ over an algebraically closed field $F$ of characteristic zero is isomorphic to $G_n(\b)$ or $G_n(\b)^{\mathrm{op}}$ for some $\b\in F$. 
\end{thm}

The set of all modules $G_n(\b)$ and $G_n(\b)^{\mathrm{op}}$ for all $\b\in F$ does not contain any pair of isomorphic modules. 
In fact, the modules  $G_n(\a)$ and $G_n(\b)^{\mathrm{op}}$ cannot be isomorphic to each other,  since the image of the identity element $1$ must be $1$. Likewise, the modules $G_n(\a)$ and $G_n(\b)$ cannot be isomorphic to each other for different values of $\a$ and $\b$, because the action of $h(1)$ is determined by these parameters.

\section{Representations of $G_n$ as a superalgebra with a Jordan bracket}.

\hspace{\parindent}

By \eqref{Jord}, we can turn every $G_n$-supermodule with a contact Lie bracket into a $G_n$-supermodule with a Jordan bracket. 
Applying this to the supermodules $G_n(\b), \b\in F$, we get the $G_n$-supermodules $\widetilde {G_n(\b)}$ with a Jordan bracket. Then $\widetilde {G_n(\b)}$ is an isomorphic  copy of $G_n$ and the actions of $G_n$ on $\widetilde {G_n}$ are defined by 
\bes
\widetilde v\cdot a&=&\widetilde{v\cdot a},\\
\la \widetilde v,a\ra&=&\widetilde{\{v,a\}}-\b \widetilde{v\cdot a},
\ees
for all $v\in V,\, a\in G_n$.   

Since $G_n(\b)$ is an irreducible $G_n$-supermodule with a contact Lie bracket it follows that 
 $\widetilde {G_n(\b)}$ is an irreducible $G_n$-supermodule with a Jordan bracket, and, consequently,  
 $\mathrm{Kan}(\widetilde {G_n(\b)})$ is an irreducible supermodule over $\mathrm{Kan}(G_n)$. 

\begin{thm}\label{tt}
Every irreducible Jordan supermodule over the superalgebra $\mathrm{Kan}(G_n)$ is isomorphic to one of the supermodules $\mathrm{Kan}(\widetilde {G_n(\b)}),\, \b\in F,$ or to their opposite supermodules.
\end{thm}
\prf
The irreducible Jordan supermodules over the superalgebra $\mathrm{Kan}(G_n)$ were classified in \cite{FSh, MZ}. 
Every irreducible Jordan supermodule over $\mathrm{Kan}(G_n)$ \cite{FSh} is isomorphic to $M_{\a}$ for some $\a\in F$ or to its opposite. Let us first recall the description of $M_{\a}$ from  \cite[Theorem 4.3]{FSh}. 

The subsets $I$ of $I_n$ are considered as ordered subsets in \cite{FSh}. For example, the subsets $\{1,2,3,4\}$ and $\{4,3,1,2\}$ are different. 
The elements $w_I$ are defined for any ordered subset $I$ of $I_n$. 
 If $\s$ is a permutation of elements of $I$ then $w_{\s(I)}=\mathrm{sgn}(\s) w_I$, where $\mathrm{sgn}(\s)$ is the sign of $\s$.  In particular, $w_{\{4,3,1,2\}}=-w_{\{1,2,3,4\}}$. 
 We say that $I\subseteq I_n$ is {\em increasing} if it does not contain any inversions. 

Let $W_{\a}$ be a vector space with a linear basis $w_I$, where $I$ runs over increasing subsets of $I_n$. Then 
\bes
M_{\a}=W_{\a}\oplus \overline{W_{\a}}, 
\ees
where $\overline{W_{\a}}$ is a copy of the vector space $W_{\a}$. 

Let $I,J$  be arbitrary ordered subsets $I_n$ such that $J=\{j_1,\ldots,j_{s_1},j_{s_1+1},\ldots,j_{s_1+s_2}\}$, $I=\{i_1,\ldots,i_{k-s_1},j_{s_1},\ldots,j_1\}$, and $|I\cap J|=s_1$. 
The action of $\mathrm{Kan}(G_n)$ on $M_{\a}$ is defined by 

$(a)$ $w_I \ e_J=\left\{\begin{array}{ccc}
0 & {\rm if} &J\not \subseteq I,\\
w_{I\setminus J}     & {\rm if} & J\subseteq I,\\
\end{array}\right.$

$(b)$ $w_I \ \overline{e_J}=\left\{\begin{array}{ccc}
0 & {\rm if} &J\not\subseteq I,\\
\overline{w_{I\setminus J}}    & {\rm if} & J\subseteq I,\\
\end{array}\right.$

$(c)$ $\overline{w_I} \ e_J=\left\{\begin{array}{ccc}
0 & {\rm if} &J\not\subseteq I,\\
(-1)^{|J|} \overline{w_{I\setminus J}}    & {\rm if} & J\subseteq I,\\
\end{array}\right.$

$(d)$ $\overline{w_I} \ \overline{e_J}=\left\{\begin{array}{ccc}
0 & {\rm if} & |J\setminus I|\geq 2,\\
(-1)^{|I\cap J|} w_{I'}    & {\rm if} & |J\setminus I|=1 (s_2=1),\\
(-1)^{|J|-1} \a (|J|-1)w_{I\setminus J}    & {\rm if} & J\subseteq I,\\
\end{array}\right.$

where $I'=\{i_1,\ldots,i_{k-s_1},j_{s_1+1}\}$.

First show that $W_{\a}$ is a $G_n$-supermodule with a Jordan bracket with respect to the actions $(a)$ and 
\bes
\langle w_I, e_J \rangle =(-1)^{|J|}\overline{w_I} \ \overline{e_J}=\left\{\begin{array}{ccc}
0 & {\rm if} & |J\setminus I|\geq 2,\\
- w_{I'}    & {\rm if} & |J\setminus I|=1,\\
- \a (|J|-1)w_{I\setminus J}    & {\rm if} & J\subseteq I.\\
\end{array}\right.
\ees

Notice that if it holds, then $(a)$--$(d)$ immediately implies that $M_{\a}=\mathrm{Kan}\,W_{\a}$ over $\mathrm{Kan}(G_n)$. 

By (\ref{Jord}), we get 
\bes
\la 1,e_J\ra=\{1,e_J\}-1/2(1\{1,e_J\}-\{1,1\}1)=0.
\ees
 Notice also that $\la w_I,1\ra=\a w_I$ by the definition of the bracket. 

In order to show that $W_{\a}$ is a $G_n$-supermodule with a Jordan bracket, by (\ref{cont}), it is sufficient to show that $W_{\a}$ is a $G_n$-supermodule with a contact Lie bracket with respect to the actions $(a)$ and 
\bes
\{w_I, e_J\} =\langle w_I, e_J\rangle+ (w_I\la 1,e_J\ra+\la w_I,1\ra e_J) \\
=(-1)^{|J|}\overline{w_I} \ \overline{e_J}+\a w_Ie_J
=\left\{\begin{array}{ccc}
0 & {\rm if} & |J\setminus I|\geq 2,\\
- w_{I'}    & {\rm if} & |J\setminus I|=1,\\
- \a (|J|-2)w_{I\setminus J}   & {\rm if} & J\subseteq I.\\
\end{array}\right.
\ees

 Let $I=\{i_1,\ldots,i_k\}$ be an increasing subset of $I_n$ and $j\in I_n$. For any ordered subset $J\subseteq I_n$ denote by $c(J)$ the increasing set obtained from $J$. If $j\in I_n\setminus I$,  then denote by $s(I,j)$ the number of inversions in the sequence $i_1,\ldots,i_k,j$ as in Lemma \ref{lem2}. 

Accurately following the definitions, one can easily check that  

(1) $w_I\cdot e_j=0$ for any $j\not\in I$,

(2) $w_I\cdot e_j=(-1)^{k-s}w_{c(I\setminus\{j\})}$ for  $j=i_s\in I$,

(3) $\{w_I,e_j\}=-(-1)^{s(I,j)}w_{c(I\cup\{j\})}$ for $j\not\in I$,

(4) $\{w_I,e_j\}=(-1)^{k-s} \a w_{c(I\setminus \{j\})}$ for $j=i_s\in I.$

Notice that there is only a sign difference between these table of actions and the table of actions given in Lemma \ref{lem2}. Obviously, 
the map $\phi: V_{\a}\to W_{\a}$ determined by $\phi(v_I)=(-1)^{|I|}w_I$ is an isomorphism $G$-supermodules with a contact Lie bracket. 
Consequently, $W_{\a}$ is isomorphic to $G_n(\a/2)$ or $G_n(\a/2)^{\mathrm{op}}$ by Lemma \ref{lem4}. 

Moreover, by the convertibility of supermodules determined by (\ref{Jord}) and (\ref{cont}), the $G_n$-supermodule $W_{\a}$ with a Jordan bracket is isomorphic to either $\widetilde{G_n(\a/2)}$ or $\widetilde{G_n(\a/2)}^{\mathrm{op}}$. Finaly, this implies that $\mathrm{Kan}(G_n(\b))$-supermodule $M_{\a}= \mathrm{Kan}(W_{a})$ is isomorphic to either $\mathrm{Kan}(\widetilde {G_n(\a/2)})$ or $\mathrm{Kan}(\widetilde {G_n(\a/2)}^{\mathrm{op}})\simeq \mathrm{Kan}(\widetilde {G_n(\a/2)})^{\mathrm{op}}$. 
\ctd

In view of Theorem \ref{tt} we formulate the following
\begin{conj}
Let $A=\la A,\cdot,\{,\}\ra$ be a simple dot-bracket superalgebra with a Jordan bracket $\{,\}$ (that is,  without ideals invariant with respect to the bracket). Then every irreducible Jordan supermodule over the superalgebra $\mathrm{Kan}(A)$ has the form $\mathrm{Kan}(V)$, where $\la V,\{,\}\ra $ is an irreducible dot-bracket  supermodule with a Jordan bracket over the superalgebra $A$.
\end{conj}

\section*{Acknowledgments}

The authors would like to thank the Shenzhen International Center for Mathematics (SUSTech) for
its hospitality and excellent working conditions, where some part of this work has been done.
The first author is supported by FAPESP  grant 2024/14914-9 and CNPq grant 305196/2024-3, and also by IMC of SUSTech.
 The second author is supported by grant AP23486782 from the Ministry of Science and Higher Education of the Republic of Kazakhstan.

 \end{document}